\documentclass{article}[12pt]

\usepackage{latexsym}
\usepackage{amsmath}
\usepackage{amssymb}

\begin{document}

\title{Two classes of explicitly solvable sextic equations}

\author{Francesco Calogero$^{a,b}$\thanks{e-mail: francesco.calogero@roma1.infn.it}
\thanks{e-mail: francesco.calogero@uniroma1.it}
 , Farrin Payandeh$^c$\thanks{e-mail: farrinpayandeh@yahoo.com}
 \thanks{e-mail: f$\_$payandeh@pnu.ac.ir}}

\maketitle   \centerline{\it $^{a}$Physics Department, University of
Rome "La Sapienza", Rome, Italy}

\maketitle   \centerline{\it $^{b}$INFN, Sezione di Roma 1}

\maketitle

\maketitle   \centerline{\it $^{c}$Department of Physics, Payame
Noor University, PO BOX 19395-3697 Tehran, Iran}

\maketitle

\begin{abstract}

The \textit{generic} monic polynomial of \textit{sixth} degree features $6$
\textit{a priori arbitrary} coefficients. We show that if these $6$
coefficients are \textit{appropriately} defined---in two different ways---in
terms of $5$ \textit{arbitrary} parameters, then the $6$ roots of the
corresponding polynomial can be \textit{explicitly} computed in terms of
\textit{radicals} of these parameters. We also report the $2$ \textit{%
constraints }on the $6$ coefficients of the polynomial implied by the fact
that they are so defined in terms of $5$ \textit{arbitrary} parameters; as
well as the \textit{explicit} determination of these $5$ parameters in terms
of the $6$ coefficients of the \textit{sextic} polynomial.

\end{abstract}

\section{Introduction}

The task of computing the roots of a given polynomial has been a fundamental
problem and a significant engine of progress in mathematics. An important
breakthrough occurred about five centuries ago, with the discovery of a
technique to find the roots of a generic polynomial of \textit{third} degree
(for a terse description of this development see, for instance, the item
"Cubic equation" in Wikipedia). The second major progress occurred about two
centuries ago and was due to Paolo Ruffini, Niels Henrick Abel and \'{E}%
variste Galois: it was the proof that the zeros of a \textit{generic}
polynomial of any degree $N>4$ \textit{cannot} be represented---in terms of
the coefficients of that polynomial---by a formula involving \textit{only
radicals:} this breakthrough opened the way to what is now called "Galois
theory" (see for instance the item "Galois Theory" in Wikipedia). This
development of course does not exclude that there exist \textit{specific}
polynomials of arbitrary degree $N$ \textit{all} roots of which can be
\textit{explicitly} computed in terms of elementary functions: for instance
consider the \textit{monic} polynomial of \textit{arbitrary} degree $6,$
\begin{subequations}
\begin{equation}
P_{N}\left( z\right) =\prod\limits_{n=1}^{N}\left( z-z_{n}\right)
=z^{N}+\sum_{m=0}^{N-1}\left( c_{m}z^{m}\right) ~,  \label{PNz}
\end{equation}%
where we now imagine the $N$ roots $z_{n}$ to be \textit{a priori
arbitrarily assigned} (hence to be \textit{known}), and the $N$ coefficients
$c_{m}$ to be then computed in terms of them---as easily implied by the
simultaneous validity of the two expressions of the polynomial (\ref{PNz})
as a product and a sum, see for instance the item "Vieta's formulas" in
Wikipedia. But of course any attempt to \textit{invert} the Vieta's formulas
in order to obtain the $N$ zeros $z_{n}$ from the $N$ coefficients $c_{m}$
would eventually require the solution of an algebraic equation of degree $N$.

Nevertheless the identification of classes of polynomials of degree $N>4$%
---defined by assigning their $N$ coefficients $c_{m},$ see (\ref{PNz}%
)---which do allow their $N$ zeros to be computed by \textit{radicals} is an
interesting mathematical topic. For instance a rather recent example---based
on Galois theory---of such an endeavour for polynomials of \textit{sixth}
degree is provided by the paper \cite{H2020} (easily reachable via Google).
An analogous endeavour---but based on more elementary mathematics---is
reported in the present paper, where \textit{two} classes of \textit{sextic
monic} polynomials are identified, which allow the computation of their $6$
zeros by \textit{radicals }from their $6$ coefficients, provided these are
defined by \textit{explicitly} provided formulas in terms of $5$ \textit{%
arbitrary} parameters. These findings are relatively trivial in the context
of Galois theory and are presumably already implied by the results reported
in \cite{H2020}; but are obtained below by more elementary means, and they
are also much simpler, as indicated by the fact that each of the \textit{%
explicit} formulas written in the present paper require only \textit{one}
line to be displayed, while several of those displayed in the paper \cite%
{H2020} takes one or more pages.

Hereafter the \textit{generic monic} polynomial of \textit{sixth} degree is
defined as follows:
\end{subequations}
\begin{subequations}
\begin{equation}
P_{6}\left( x\right) =x^{6}+\sum_{n=0}^{5}\left( c_{n}x^{n}\right) ~,
\label{Pol6}
\end{equation}%
with the $6$ coefficients $c_{n},$ $n=0,1,2,3,4,5$ \textit{a priori
arbitrary }(except for the conditions mentioned below); and $z_{\lambda \mu
} $ are the $6$ roots of this polynomial,
\begin{equation}
P_{6}\left( z_{\lambda \mu }\right) =0~,~~~\lambda =1,2,~~\mu =1,2,3~.
\label{zeros}
\end{equation}

\section{Results: first model}

\textbf{Proposition 2-1}. Assume that the $6$ coefficients $c_{n}$ of the
\textit{sextic} polynomial (\ref{Pol6}) may be expressed as follows in terms
of the $5$ \textit{arbitrary} parameters $a_{0},$ $a_{1},$ $a_{2},$ $b_{0},$
$b_{1}$:
\end{subequations}
\begin{subequations}
\label{cc}
\begin{equation}
c_{5}=2a_{2}~,  \label{cc5}
\end{equation}%
\begin{equation}
c_{4}=2a_{1}+\left( a_{2}\right) ^{2}~,  \label{cc4}
\end{equation}%
\begin{equation}
c_{3}=2a_{0}+2a_{1}a_{2}~+b_{1}~,  \label{cc3}
\end{equation}%
\begin{equation}
c_{2}=\left( a_{1}\right) ^{2}+\left( 2a_{0}+b_{1}\right) a_{2}~,
\label{cc2}
\end{equation}%
\begin{equation}
c_{1}=\left( 2a_{0}+b_{1}\right) a_{1}~,  \label{cc1}
\end{equation}%
\begin{equation}
c_{0}=\left( a_{0}\right) ^{2}+a_{0}b_{1}+b_{0}~.  \label{cc0}
\end{equation}%
Then the $6$ roots $z_{\lambda \mu }$ of the \textit{sextic} polynomial (\ref%
{Pol6}) are \textit{explicitly} given, in terms of the $5$ parameters $%
a_{0}, $ $a_{1},$ $a_{2},$ $b_{0},$ $b_{1},$ by the following definitions:
the $6$ numbers $z_{\lambda \mu }$ are the $3$ roots (with $\lambda =1,2$
and $\mu =1,2,3$) of the following $2$ \textit{cubic} equation in $z$,
\end{subequations}
\begin{subequations}
\label{Eqzy}
\begin{equation}
z^{3}+a_{2}z^{2}+a_{1}z+a_{0}=y_{\lambda }~,~~~\lambda =1,2~,  \label{Equaz}
\end{equation}%
where $y_{\lambda }$ is one of the $2$ roots of the following \textit{%
quadratic }equation in $y$,
\begin{equation}
y^{2}+b_{1}y+b_{0}=0~.~\blacksquare  \label{Equay}
\end{equation}

\textbf{Remark 2-1}. The \textit{quadratic} respectively \textit{cubic}
equations (\ref{Equay}) respectively (\ref{Equaz}) can of course be solved
\textit{explicitly}:
\end{subequations}
\begin{equation}
y_{\lambda }=\left\{ -b_{1}+\left( -1\right) ^{\lambda }\sqrt{\left(
b_{1}\right) ^{2}-4b_{0}}\right\} /2~,~~~\lambda =1,2~,
\end{equation}%
while the $3$ roots $z_{\lambda \mu }$ of the \textit{cubic} equation (\ref%
{Equaz}) are given by the well-known "Cardano" formulas (see again the item
"Cubic equation" in Wikipedia). The resulting \textit{explicit} formula
expressing the $6$ zeros $z_{\lambda \mu }$ in terms of the $5$ parameters $%
a_{0},$ $a_{1},$ $a_{2},$ $b_{0},$ $b_{1}$ involves---in a nested way---only
\textit{square} and \textit{cubic} roots; it is of course a bit more
complicated than the Cardano formulas, as the reader who takes the
trouble---indeed, an easy task---to write it out shall easily find out. $%
\blacksquare $

But, as indicated above, the task of obtaining---from the \textit{assignment}
of a number of parameters---a corresponding set of both the \textit{%
coefficients} of a polynomial and its \textit{zeros} may be a relatively
easy task. Less trivial is the task to assign \textit{a priori} the $N$
\textit{coefficients} $c_{n}$ of a monic polynomial of degree $N$ (see (\ref%
{PNz})) and to then find a (\textit{generally smaller}) number of \textit{%
parameters} which determine---as it were, \textit{a posteriori}---via
\textit{explicit} formulas both these preassigned $N$ coefficients $c_{n}$
and the $N$ \textit{zeros} of the corresponding monic polynomial, as well as
the \textit{explicit} formulas displaying the corresponding \textit{%
constraints} implied by these assignments on the $N$ coefficients $c_{n}$.

The following proposition provides---for $N=6$---such findings, which
complement those reported in \textbf{Proposition 2-1}.

\textbf{Proposition 2-2}. If the $6$ parameters $c_{n}$ are expressed in
terms of the $5$ parameters $a_{0},$ $a_{1},$ $a_{2},$ $b_{0},$ $b_{1}$ by
the $6$ formulas (\ref{cc}), then the $5$ parameters $a_{0},$ $a_{1},$ $%
a_{2},$ $b_{0},$ $b_{1}$ are themselves expressed as follows in terms of the
$6$ coefficients $c_{n}$:
\begin{subequations}
\label{ab}
\begin{equation}
a_{2}=\left( c_{5}/2\right) ~,~~~  \label{a2}
\end{equation}%
\begin{equation}
a_{1}=\left[ c_{4}-\left( a_{2}\right) ^{2}\right] /2=\left[ 4c_{4}-\left(
c_{5}\right) ^{2}\right] /8~,  \label{a1}
\end{equation}%
\begin{equation}
2a_{0}+b_{1}=c_{3}-2a_{1}a_{2}=\left[ c_{2}-\left( a_{1}\right) ^{2}\right]
/a_{2}=c_{1}/a_{1}~,  \label{a0b1}
\end{equation}%
\begin{equation}
b_{0}=c_{0}-\left( a_{0}+b_{1}\right) a_{0}~;  \label{b0}
\end{equation}%
and the $6$ coefficients $c_{n}$ satisfy the following $2$ \textit{%
constraints}:
\end{subequations}
\begin{subequations}
\label{Con}
\begin{equation}
c_{1}=a_{1}\left( 2a_{0}+b_{1}\right) =\left[ 4c_{4}-\left( c_{5}\right) ^{2}%
\right] \left\{ c_{3}-\left[ 4c_{4}-\left( c_{5}\right) ^{2}\right]
c_{5}/8\right\} /8~,  \label{Con1}
\end{equation}%
\begin{equation}
c_{2}=a_{2}c_{3}+a_{1}\left[ a_{1}-2\left( a_{2}\right) ^{2}\right]
=c_{3}c_{5}/2+\left[ 4c_{4}-\left( c_{5}\right) ^{2}\right] \left[
4c_{4}-5\left( c_{5}\right) ^{2}\right] /64~.  \label{Con2}
\end{equation}%
Note that the $2$ parameters $a_{2}$ and $a_{1}$ are given \textit{explicitly%
} in terms of the coefficients $c_{4}$ and $c_{5}$ by the eqs. (\ref{a2})
and (\ref{a1}); while the $3$ equalities (\ref{a0b1}) imply (via (\ref{a2})
and (\ref{a1})) the $2$ \textit{constraints} (\ref{Con}) on the $6$
coefficients $c_{n}$---which express \textit{explicitly} the $2$
coefficients $c_{1}$ and $c_{2}$ in terms of the $3$ coefficients $c_{3},$ $%
c_{4},$ $c_{5}$. Then---once these $2$ \textit{constraints} are
satisfied---the $3$ equalities (\ref{a0b1}), together with eq. (\ref{b0}),
provide the explicit determination of the parameters $b_{0}$ and $b_{1}$ in
terms of the parameters $c_{n}$ and of the parameter $a_{0}$---this easy
task amounts to solving a system of $2$ \textit{linear} equations---while $%
a_{0}$ remains as a \textit{free} parameter; this freedom might be used to
simplify all the above formulas, for instance by assuming that $a_{0}$
vanishes or that $a_{0}=-b_{1}/2$ implying $c_{1}=0$ (see (\ref{a0b1})), but
of course at the cost of decreasing the generality of these findings. $%
\blacksquare $~

Proofs of \textbf{Propositions 2-1} and \textbf{2-2 }are provided in the
\textbf{Appendix}.

\bigskip

\section{Results: second model}

The findings reported in this \textbf{Section 3} are \textit{analogous}, but
\textit{different}, from those reported in \textbf{Section 2}; accordingly,
the variables and parameters used in this \textbf{Section 3 }are \textit{%
different} from those having the \textit{same name} in \textbf{Section 2},
although they play \textit{analogous} roles.

We only report below these new findings, without detailing their derivation;
which is quite \textit{analogous} to that described above and below (see
\textbf{Section 2} and\textbf{\ Appendix A}), and may be recommended as an
interesting exercise for the enterprising reader (clue: compare the eqs. (%
\ref{Eqzy}) to the eqs. (\ref{Eqyz}), see below).

\textbf{Proposition 3-1}. Assume that the $5$ coefficients $c_{n}$ may be
expressed as follows in terms of the $5$ \textit{arbitrary} parameters $%
a_{0},$ $a_{1},$ $a_{2},$ $b_{0},$ $b_{1}$:
\end{subequations}
\begin{subequations}
\label{ccc}
\begin{equation}
c_{5}=3b_{1}~,  \label{ccc5}
\end{equation}%
\begin{equation}
c_{4}=a_{2}+3\left[ b_{0}+\left( b_{1}\right) ^{2}\right] ~,  \label{ccc4}
\end{equation}%
\begin{equation}
c_{3}=\left[ 2a_{2}+6b_{0}+\left( b_{1}\right) ^{2}\right] b_{1}~,
\label{ccc3}
\end{equation}%
\begin{equation}
c_{2}=a_{1}+a_{2}\left[ 2b_{0}+\left( b_{1}\right) ^{2}\right] +3b_{0}\left[
b_{0}+\left( b_{1}\right) ^{2}\right] ~,  \label{ccc2}
\end{equation}%
\begin{equation}
c_{1}=a_{1}b_{1}+2a_{2}b_{0}b_{1}+3\left( b_{0}\right) ^{2}b_{1}~,
\label{ccc1}
\end{equation}%
\begin{equation}
c_{0}=a_{0}+a_{1}b_{0}+a_{2}\left( b_{0}\right) ^{2}+\left( b_{0}\right)
^{3}~.  \label{ccc0}
\end{equation}%
Then the $6$ roots $z_{\lambda \mu }$ (with $\lambda =1,2$ and $\mu =1,2,3$)
of the \textit{sextic} polynomial (\ref{Pol6}) are \textit{explicitly}
given, in terms of the $5$ parameters $a_{0},$ $a_{1},$ $a_{2},$ $b_{0},$ $%
b_{1},$ by the following definitions: $z_{\lambda \mu }$ is one of the $2$
roots of the following $3$ \textit{quadratic }equations in $z$,
\end{subequations}
\begin{subequations}
\label{Eqyz}
\begin{equation}
z^{2}+b_{1}z+b_{0}=y_{\mu }~,~~~\mu =1,2,3~,
\end{equation}%
where $y_{\mu }\ $is one of the $3$ roots of the following \textit{cubic} in
$y$:
\begin{equation}
y^{3}+a_{2}y^{2}+a_{1}y+a_{0}=0~.~\blacksquare
\end{equation}

\textbf{Proposition 3-2}. If the $6$ parameters $c_{n}$ are expressed in
terms of the $5$ parameters $a_{0},$ $a_{1},$ $a_{2},$ $b_{0},$ $b_{1}$ by
the $6$ formulas (\ref{ccc}), then the $5$ parameters $a_{0},$ $a_{1},$ $%
a_{2},$ $b_{0},$ $b_{1}$ are themselves expressed as follows in terms of the
$6$ coefficients $c_{n}$:
\end{subequations}
\begin{subequations}
\label{32ab}
\begin{equation}
b_{1}=c_{5}/3~,  \label{32b1}
\end{equation}%
\begin{equation}
a_{2}=c_{4}-\left( c_{5}\right) ^{2}/3-3b_{0}  \label{32a2}
\end{equation}%
\begin{equation}
a_{1}=c_{2}-2\left[ 3c_{4}-9b_{0}-\left( c_{5}\right) ^{2}\right] \left[
18b_{0}+\left( c_{5}\right) ^{2}\right] /27-b_{0}\left[ 9b_{0}+\left(
c_{5}\right) ^{2}\right] /3~,  \label{32a1}
\end{equation}%
\begin{equation}
a_{0}=c_{0}-a_{1}b_{0}-a_{2}\left( b_{0}\right) ^{2}-\left( b_{0}\right)
^{3}~;  \label{32a0}
\end{equation}%
while the $6$ coefficients $c_{n}$ are required to satisfy the following $2$
\textit{constraints}:
\end{subequations}
\begin{subequations}
\label{32Con}
\begin{equation}
27c_{3}-18c_{4}c_{5}+5\left( c_{5}\right) ^{3}=0~,  \label{32Con1}
\end{equation}%
\begin{equation}
c_{1}=\left[ 27c_{2}-3c_{4}\left( c_{5}\right) ^{2}+\left( c_{5}\right) ^{4}%
\right] c_{5}/81~.  \label{32Con2}
\end{equation}%
Note that in this case we wrote out \textit{explicitly}---in terms of the
parameters $c_{n}$---the expressions of the $3$ parameters $b_{1},$ $a_{2}$
and $a_{1}$, and that an \textit{explicit} expression of the \ parameter $%
a_{0}$ is also implied by (\ref{32a0}) via (\ref{32a1}) and (\ref{32a2});
while the parameter $b_{0}$ remains as a \textit{free} parameter (the same
role played by the parameter $a_{0}$ in \textbf{Proposition 2-2}). $%
\blacksquare $

\bigskip

\section{Outlook}

Obvious generalizations of the approach employed in this paper may be used
to get analogous results for polynomials of degree $N>6,$ especially
whenever $N=2^{p_{1}}3^{p_{2}}$ with $p_{1}$ and $p_{2}$ \textit{arbitrary
nonnegative integers }(for instance for $N=8$ or $N=9$); but then the number
of restrictions on the coefficients of these polynomials for the
applicability of this approach shall of course grow as $N$ grows.

\bigskip

\section{Appendix}

In this Appendix we prove the results reported in \textbf{Section 2}.

To derive the results reported in \textbf{Proposition 2-1} all one needs to
do is to replace $y$ in eq. (\ref{Equay}) by the expression $%
z^{3}+b_{2}z^{2}+b_{1}z+b_{0}$ (see (\ref{Equaz})), expand the resulting
expression in powers of $z$, and identify the coefficients of the resulting
\textit{sextic} equation in $z$ with the coefficients $c_{n}$,$\ $see (\ref%
{Pol6}).

The results reported in \textbf{Proposition 2-2} are easy consequences of
the formulas (\ref{cc}). It is plain that the $3$ eqs. (\ref{a2}), (\ref{a1}%
) respectively (\ref{b0}) are implied by the $3$ eqs. (\ref{cc5}), (\ref{cc4}%
) respectively (\ref{cc0}). Next, it is easily seen that the $3$ equalities (%
\ref{a0b1}) are implied by the $3$ eqs. (\ref{cc3}), (\ref{cc2}) and (\ref%
{cc1}). And the fact that these $3$ equalities imply the $2$ \textit{%
constraints} (\ref{Con}) is then obvious, as well as the fact that they
allow the \textit{explicit} computation of the parameters $b_{0}$ and $b_{1}$
while leaving $a_{0}$ as a free (undetermined) parameter; as mentioned in
\textbf{Proposition 2-2}, which is thereby proven.

\end{subequations}


\bigskip

\begin{thebibliography}{9}
\bibitem{H2020} T. R. Hagedorn, "General formulas for solving solvable
sextic equations", J. Algebra \textbf{233}, 704-757 (2000).
\end{thebibliography}
\end{document}